\documentclass{amsart}
\usepackage{amssymb,amsmath,amsthm,longtable}
\usepackage[dvips]{graphicx}

\newcommand{\D}{{\mathbb D}}
\newcommand{\C}{{\mathbb C}}
\newcommand{\E}{{\mathbb E}}

\newcommand{\R}{{\mathbb R}}

\def\br#1{\left(#1\right)}
\def\brb#1{\left[#1\right]}
\def\brs#1{\left\{#1\right\}}
\newcommand{\K}{K}

\theoremstyle{remark}

\theoremstyle{remark}

\theoremstyle{definition}

\theoremstyle{plain}

\date{}
\title{Integral means spectrum of random conformal snowflakes}
\author{D. Beliaev }
\begin{document}
\maketitle
%\begin{abstract}

%\end{abstract}

\section{Introduction}

It is known that extremal configurations in many important 
problems in classical complex analysis exhibit complicated fractal structure. 
This makes such problems extremely difficult. 
The classical example is the coefficient problem 
for the class of bounded univalent functions.
Let $\phi(z)=z+a_2 z^2+a_3 z^3+\dots$ be a bounded univalent map in the unit disc. 
One can ask what are the maximal possible values of coefficients $a_n$, 
especially when $n$ tends to infinity. We define $\gamma$ as the best 
possible constant such that $|a_n|$ decays as $n^{\gamma-1}$. In \cite{CaJo} 
Carleson and Jones showed that this problem is related to 
another classical problem about the growth rate of the length of Greens' lines. 
In particular, they showed that the extremal configurations for both of this problems 
should be of a fractal nature.

During the last decade it became clear that the right language for
these problems, as well as many other classical problems, 
is {\em the maltifractal analysis}. It turned out that all these problems
could be reduced to the problem of finding the maximal value of 
{\em the integral means spectrum}.

In the recent paper \cite{BeSmsnow} S.~Smirnov and the author introduced 
 and studied a new class of random fractals,
the so-called {\em random conformal snowflakes}. In particular, they
proved the {\em fractal approximation} for this class, which
 means that one can find conformal snowflakes with
spectra arbitrary close to the maximal possible spectrum. In this
paper we report on our search for snowflakes with large spectrum.

The paper organized as follows: in the introduction we give some
basic information about integral means spectrum, define random conformal
snowflakes, and state the main facts about spectrum of the snowflakes. 
In the Section \ref{sec:num} we give numerical estimates
of the spectra of several snowflakes for different values of
parameters. In the last Section we give rigorous lower bound for the
spectrum at $t=1$.

\subsection{Integral means spectrum}

Here we briefly sketch the necessary definitions and the background. For the
detailed discussion of the subject we recommend surveys
\cite{Makarov,Pommerenke97,BeSmECM} and books
\cite{Pommerenke75,Pommerenke92}.

Let $\Omega\in \hat \C$ be a simply connected domain in the complex
sphere which contains infinity. By the Riemann uniformization
theorem, there is a conformal map $\phi$ from the complement of the
unit disc onto $\Omega$ such that $\phi(\infty)=\infty$.  The {\em
integral means spectrum} is defined as
$$
\beta(t)=\beta_\Omega(t)=\beta_\phi(t)=\limsup_{r\to1+}
\frac{\log\int_{0}^{2\pi}|\phi'(re^{i\theta})|^t
d\theta}{|\log(r-1)|},~t\in\R~.
$$
The  {\em universal integral means spectrum} is
$$
B(t)=\sup \beta(t),
$$
where supremum is over all simply connected domains $\Omega$.

The universal spectrum plays the central role in the Geometric
Function Theory. A lot of major problems in the field can be stated in terms of $B(t)$.
In particular, Brennan's conjecture \cite{Brennan} asserting that $|\psi'|^{4-\epsilon}$ is
integrable for any $\epsilon>0$ and any conformal $\psi:\Omega\to\D$ 
is equivalent to $B(-2)=1$. The
Carleson-Jones conjecture, $\gamma=1/4$, is equivalent to $B(1)=1/4$. 
In 1996 Kraetzer formulated
the ultimate conjecture on integral means:
\begin{eqnarray*}
B(t)&=& t^2/4, \quad |t|\le 2,\\
B(t)&=& |t|-1, \quad |t|>2.
\end{eqnarray*}
This conjecture is based on the above mentioned conjectures,
computer experiments and the believe that the spectrum should be
relatively simple function. Kraetzer performed computer experiments with
quadratic Julia sets. The results were within $5$ percents
from $t^2/4$. We should point out that the experiments were completely non-rigorous, and, as far as we know, the error can not be estimated using available techniques.

In the present paper we give a computer
assisted lower bound on the universal spectrum which are  within
$4-7$ percents from $t^2/4$ but they are ``semi-rigorous''. By this
we mean that we use computer only for numerical integration. For
numerical integration of an explicit function one can write
estimates of the error term and this will give a rigorous estimate.
We give an estimate of the error term only for $t=1$ (this case is
of special interest and estimates are a bit simpler). We prove that 
$B(1)>0.23$ which is a significant improvement over previously known
$B(1)>0.17$. 

Essentially the same  argument as for the case $t=1$ could be used to
estimate the error for other values of $t$, although we are not doing it here.

\subsection{Conformal snowflake}
\label{sec:def}

 We denote  by
$\Sigma'$ the class of univalent maps $\phi$ from the
complement of the unit disc $\D_-$ into itself such that
$\phi(\infty)=\infty$ and $\phi'(\infty)\in \R$. Let $K_n\phi$ be
the Koebe transform of $\phi\in \Sigma'$:
$$
(\K_n\phi)(z)=(\K_n\phi)(z)=\sqrt[n]{\phi(z^n)}.
$$
 We denote the conjugation by rotation by $\phi_\theta(z)=e^{i\theta}\phi(z e^{-i\theta})$.

To construct a random snowflake we need two components: a building
block $\phi \in \Sigma'$ and an integer $k\ge 2$. The $n$-th
approximation to the snowflake is
\begin{eqnarray*} f_{n}(z)=f_{n-1}(K_{k^n}\phi_{\theta_n}(z))=
\phi_{\theta_0}(\phi_{\theta_1}^{1/k}(\dots\phi_{\theta_n}^{1/k}(z^{k^n})\dots),
\end{eqnarray*}
where $\theta_j$ are independent variables uniformly distributed
between $0$ and $2\pi$. The conformal snowflake $f$ is the limit of
$f_n$. For the proof that conformal snowflakes are well defined see
\cite{BeSmsnow}.  $S=\C\setminus f(\D_-)$ and
$g=f^{-1}$ are also referred to as snowflakes.

\begin{figure}
    \centering
        \includegraphics{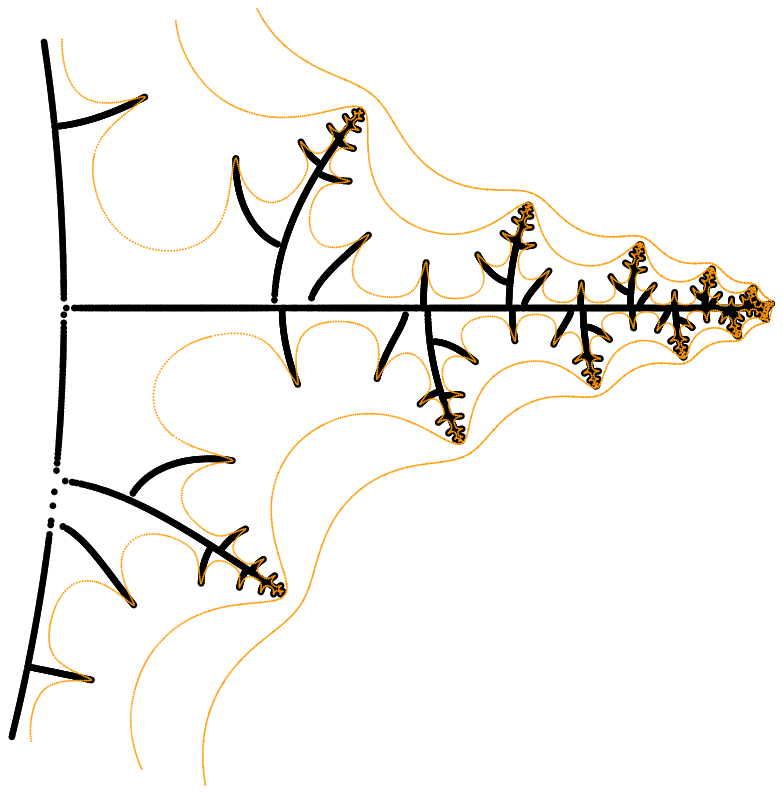}
    \caption{The image of a small boundary arc under $ f_3$ with three Green's lines.}
    \label{pic3}
\end{figure}

It turns out that the spectrum of a snowflake for a particular value
of $t$ is closely related to the spectrum of an integral operator
$P$ whose kernel is defined in terms of $k$, $t$, and $\phi$. This
operator is defined as
\begin{equation}
P\nu(r):={r^{1-\frac{(k-1)t}{k}}} \int_0^{2\pi}
\frac{\nu(\phi(r^{1/k}e^{i\theta}))}{|\phi(r^{1/k}e^{i\theta})|}
|\phi'(r^{1/k}e^{i\theta})|^{t}\frac{d\theta}{2\pi}.
\end{equation}
Let $R$ be any radius such that $D_R \subset \psi^k(D_R)$ where
$D_R=\brs{z:1<|z|<R}$ and $\psi=\phi^{-1}$.

One can show (see \cite{BeSmsnow} for the proof) that if $\lambda$
is the maximal eigenvalue of $P$ then
\begin{equation}
\label{eq:lambda}
 \beta(t)\ge \log \lambda/\log k,
\end{equation}
in particular
\begin{equation}
\label{eq:nu}
 \beta(t)\ge \min_{1\le r \le
R}\log\left(\frac{P\nu(r)}{\nu(r)}\right)/\log k
\end{equation}
for any positive test function $\nu$.

\subsection{Random fractals}

One of the main problems in the computation of the integral means
spectrum (or other multifractal spectra) is the fact that the
derivative of a Riemann map for a fractal domain depends on the
argument in a very non regular way: $\phi'$ is a ``fractal'' object
in itself. The solution to this problem is to study random fractals
for which the boundary behavior of $\phi'$ is statistically
independent of the argument. In this case it is natural to consider
the {\em average integral means spectrum:}
\begin{eqnarray*}
\bar\beta(t)&=&\sup\brs{\beta: \int_1(r-1)^{\beta-1}\int_0^{2\pi}
\E\brb{|f'(r e^{i\theta})|^t}d \theta d r=\infty}
\\
&=&\inf\brs{\beta: \int_1(r-1)^{\beta-1}\int_0^{2\pi}
\E\brb{|f'(r e^{i\theta})|^t}d \theta d r<\infty}.
\end{eqnarray*}
The average spectrum does not have to be related to the spectra of a
particular realization. We want to point out that even if $\phi$ has
the same spectrum a.s. it does not guarantee that $\bar\beta(t)$ is
equal to the a.s. value of $\beta(t)$. Moreover, it can happen that
$\bar\beta$ is not a spectrum of {\em any} particular domain. On the other hand,
Makarov's fractal approximation (\cite{Makarov}) implies that $\bar \beta(t)\le
B(t)$. Thus lower bound for any $\bar\beta$ is a lower bound for
$B$. In the remainder of the paper we will work only with average integral means
spectrum. We refer to it as ``the spectrum'' and denote it by $\beta$.

\section{Numerical estimates of the spectrum}
\label{sec:num}

In this section we give a numerical support of Kraetzer's
conjecture. To achieve this we estimate the spectrum of snowflakes with a
straight slit building blocks (see Figure 1).  We use the following scheme: first
we compute the disctretized operator $P$ and find its main
eigenvalue and eigenvector. We approximate the main eigenvector
by a nice function $\nu$ and use it as a test function in
(\ref{eq:nu}).

First we define the basic slit function
\begin{equation}
\label{slit}
\phi(z,l)=\phi_l(z)=\mu_2\br{\frac{\sqrt{\mu_1^2(z s)+l^2/(4k+4)}}{\sqrt{1+l^2/(4l+4)}}},
\end{equation}
where $s$ is a constant close to $1$, $\mu_1$ and $\mu_2$ are the
M\"obius transformation that maps $\D_-$ onto the right half plane
and its inverse:
\begin{eqnarray*}
\mu_1(z)=\frac{z-1}{z+1},\\
\mu_2(z)=\frac{z+1}{z-1}.
\end{eqnarray*}
We also need the inverse function
\begin{equation}
\psi(z,l)=\psi_l(z)=\phi(z,l)^{-1}.
\end{equation}
 The image $\phi_l(\D_-)$ is
$\D_-$ with a horizontal slit originating at $1$. The length of the
slit is $l$. The derivative of a slit map has singularities at
points mapped to $1$. But for $s>1$ these
singularities are not in $\D_-$.

%We set $s=1.002$.

Given $k$, $l$ (and $s$) we can numerically find the critical radius
$R$ and define the operator $P$. For the first step, we use the
following discrete operator which approximates $P$. Choose
sufficiently large $N$ and $M$. Let $r_n=1+(R-1)n/N$ and $\theta_m=
2\pi m /M$.  $P$ is approximated by an $N\times N$ matrix with
elements
$$
P_{n,n'}=\sum_{m}
r_n^{1-t(k-1)/k}\frac{|\phi'(r_n^{1/k}e^{i\theta_{m}})|^t}
{|\phi(r_n^{1/k}e^{i\theta_{m}})|M},
$$
where summation is over all indexes $m$ such that $r_{n'}$ is the
closest point to $|\phi(r_n^{1/k}e^{i\theta_m})|$. This defines the
discretized operator $P_N$. Let $\lambda_N$ and $V_N$ be the main
eigenvalue and the corresponding eigenvector.

It is highly plausible that the operators $P_N$ converge to $P$ and  the $\lambda_N$
converge to the main eigenvalue of $P$, but it seems to be hard to prove.
Moreover, the rate of convergence is difficult to estimate. But this crude
approximation gives us the fast test on whether the pair $\phi$ and $k$
defines a snowflake with large spectrum or not.

\begin{table}
\label{tb:beta}
\centering
\begin{tabular}{|c|c|c|c|c|c|c|}
 \hline
 \quad $t$ \qquad & \quad $k$ \quad & \quad $l$ \quad
&\quad $\log_{k} \lambda$ \qquad & \quad $\beta(t)$ \qquad &
Kraetzer &\quad $t^2/4$\qquad \\ \hline
%
%$$ & $$ &  $$ & $$ & $$ & $$ & $$ \\ \hline
%
$-2.0$ & $34$ &  $1$ & $1.262$ & $-$ & $-$ & $1.0$ \\ \hline
$-1.8$ & $34$ &  $1$ & $1.068$ & $-$ & $-$ & $0.81$ \\ \hline
$-1.6$ & $34$ &  $1$ & $0.8761$ & $-$ & $-$ & $0.64$ \\ \hline
$-1.4$ & $34$ &  $1$ & $0.6879$ & $-$ & $0.476$ & $0.49$ \\ \hline
$-1.2$ & $34$ &  $1$ & $0.5059$ & $-$ & $0.340$ & $0.36$ \\ \hline
$-1.0$ & $34$ &  $1$ & $0.3354$ & $-$ & $0.231$ & $0.25$ \\ \hline
$-0.8$ & $34$ &  $1$ & $0.1865$ & $-$ & $0.149$ & $0.16$ \\ \hline
$-0.6$ & $24$ &  $21$ & $0.0848$ & $0.0710$ & $0.085$ & $0.09$ \\
\hline
$-0.4$ & $20$ &  $25$ & $0.0377$ & $0.0352$ & $0.037$ & $0.04$ \\
\hline
$-0.2$ & $31$ &  $44$ & $0.0093$ & $0.0083$ & $0.0095$ & $0.01$ \\
\hline
$0.2$ & $5$ &  $7$ & $0.0091$ & $0.00897$ & $0.0094$ & $0.01$ \\
\hline
\ $0.4$ & $11$ &  $30$ & $0.0376 $ & $0.03767$ & $0.037$ & $0.04$ \\
\hline
$0.6$ & $14$ &  $68$ & $0.0851$ & $0.08442$ & $0.086$ & $0.09$ \\
\hline
$0.8$ & $12$ &  $67$ & $0.1514$ & $0.1511$ & $0.154$ & $0.16$ \\
\hline
$1.0$ & $13$ &  $73$ & $0.2362$ & $0.2340$ & $0.242$ & $0.25$ \\
\hline
$1.2$ & $10$ &  $67$ & $0.3425$ & $0.3350$ & $0.346$ & $0.36$ \\
\hline
$1.4$ & $8$ &  $55$ & $0.4680$ & $0.4586$ & $0.476$ & $0.49$ \\
\hline
$1.6$ & $6$ &  $39$ & $0.6137$ & $0.6091$ & $-$ & $0.64$ \\
\hline
$1.8$ & $6$ &  $39$ & $0.7790$ & $0.7713$ & $-$ & $0.81$ \\ \hline
$2.0$ & $4$ &  $21$ & $0.9548$ & $0.9296$ & $-$ & $1.0$ \\ \hline
%
%$$ & $$ &  $$ & $$ & $$ & $$ & $$ \\ \hline
\end{tabular}
%\end{center}
\caption{Spectra of nearly optimal snowflakes for different values
of $t$. Parameters $k$ and $l$ define the snowflake ($s$=1),
$\lambda$  is the main eigenvalue of $P_N$ (computed with $N=2000$
and $M=1000$), $\beta(t)$ is the lower bound given by (\ref{eq:nu}).
In the last two columns we give lower bounds obtained by Kraetzer
and values of the conjectured universal spectrum.}
\end{table}

For each value of $t$ we use this method to search for the optimal
values of $k$ and $l$. This gives us snowflakes that a supposed to
have large spectrum. To estimate the spectrum  we use the following
trick. We approximate the eigenvector $V_N$ by a simple function (or
just linearly interpolate it). The result is used as a test function
$\nu$. Than we use numerical integration to calculate $P\nu$. By
(\ref{eq:nu}) it gives us the lower bound of $\beta$. For general
values of $t$ we do it without estimates of the error term, we just
use Euler quadrature formula doubling number of nodes until the
difference is small. We also can argue that since functions seems to
be quite smooth (there is no oscillation), the precision of
numerical integration should be much better than what follows from the
standard estimates of the error term. The results of our
computations are given in the Table 1.%\ref{tb:beta}.

We want to make a few comments about negative values of $t$. The
derivative of $\phi$ has a zero of the first order at $z=1$, this
means that for $\rho=1$ the integral in the definition of $P$
diverges at $r=1$. This is the reason why we can not find a good
test function and the estimated $\lambda$ is probably far from 
the correct one for $t\le -0.8$. The
problem with $\lambda$ could be solved by increasing $N$ and $M$.
For example if we set $N=3000$ and $M=2000$ then for $t=-0.6,
\dots,-1.2$ the logarithm of $\lambda$ will be equal to $0.0847, \
0.1579,\ 0.2926$, and $\ 0.4509$ correspondingly. This means that
the value for $t=-0.6$ is probably correct even for smaller values
of $N$ and $M$ and the value for $t=-1.2$ is still far from its true
value.

\section{Spectrum at $t=1$}
\label{sec:t1}

 In this section we estimate the error term in the
numerical integration of $P\nu$. This gives a rigorous proof that
$B(1)>0.23$ which is a significant improvement over previously known
$B(1)>0.17$. We also would like to point out that this computation
is not specific for $t=1$ (though it is a bit easier in this case)
and can be done for other values as well.

For $t=1$ we study the snowflake with $l=73$ and $k=13$. The Figure
\ref{pic3} shows the image of a small arc under $ f_3$ and three
Green's lines. This snowflake has a spectrum close to $0.234$.
Unfortunately there is a singularity in the kernel of $P$ (it is
integrable, so one can still estimate the error term). This
singularity is the reason why we introduced the parameter $s$. If we
choose $s$ close to $1$ then snowflake is close to the snowflake
with $s=1$, but derivative of $\phi$ is continuous up to the
boundary. In our case we use $s=1.002$.

First we have to find the critical radius $R$ such that $D_R\subset
\psi^k(D_R)$. By symmetry of $\phi$, the critical radius is the only
positive solution of
$$
\psi^k(x)=x.
$$
This equation can not be solved explicitly, but we can solve it
numerically (we don't care about error term since we can use any upper bound for $R$). The approximate value of $R$ is $ 76.1568$.
To be on the safe side we fix $R=76.2$.

To find the test function we approximate $P$ by $P_N$ where $N=1000$
and $M=500$. The logarithm of the first eigenvalue is $0.2321$ (it
is $0.23492$ if we take $s=1$). We scale the coordinates of the main
eigenvector of $P_N$ from $[1,1000]$ to the interval $[1,R]$ and
approximate by a rational function $\nu$ which we will use as a test
function:
$$
\begin{aligned}
\nu(x)=(7.1479+8.9280 x - 0.07765 x^2+ 1.733 \times 10^{-3} x^3 -
\\ 2.0598 \times 10^{-5} x^4 +
9.5353 \times 10^{-8}x^5)/( 2.7154+ 13.2845  x).
\end{aligned}
$$

To estimate $\beta(1)$
we have to integrate $\nu(|\phi|)|\phi'|/|\phi|$. It is easy to see
that the main contribution to the derivative is given by a factor $|\phi'|/|\phi|$.
Assume for a while that $s=1$. The fraction
$|\phi'(z)/\phi(z)|$ can be written as
\begin{equation}
\label{fraction}
\frac{|z-1|}{|z|\sqrt{|(z-z_1)(z-z_2)}},
\end{equation}

where
$$
\begin{aligned}
z_1=\frac{-5033-292 i \sqrt{74}}{5625}\approx -0.894756 - 0.446556 i,
\\
z_2=\frac{-5033+292 i \sqrt{74}}{5625}\approx -0.894756 + 0.446556 i.
\end{aligned}
$$

Singular points $z_1$ and $z_2$ are mapped to $1$ and $\phi'$ has a square root type
singularity at these points. They will play essential role in all further calculations.
We introduce notation $z_1=x+i y$ and $z_2=x-i y$,
for $z$ we will use polar coordinates $z=r e^{i\theta}$.

We compute the integral of $f=|\nu(\phi)\phi'/\phi|$ using the Euler
quadrature formula based on the trapezoid quadrature  formula
$$
\int_0^{2\pi} f(x)dx \approx S_\epsilon^n(f)=
S_\epsilon(f) -\sum_{k=1}^{n-1} \gamma_{2k}\epsilon^{2k}\br{f^{(2k-1)}(2\pi)-f^{(2k-1)}(0)},
$$
where $S_\epsilon(f)$ is a trapezoid  quadrature formula with step $\epsilon$ and
$\gamma_k=B_k/k!$ were $B_k$ is the Bernoulli number.
The error term in the Euler formula is
\begin{equation}
-\gamma_{2n}\max f^{(2n)}\epsilon^{2n}2\pi.
\label{error1}
\end{equation}
In our case function $f$ is periodic and terms with higher derivatives vanish. This means
that we can use (\ref{error1}) for any $n$ as an estimate of the error
in the trapezoid quadrature formula.

Function $\phi$ has two singular points: $z_1$ and $z_2$. Derivative of $\phi$ blows up near
these points. This is why we  introduce scaling factor $s$. We can write a power series of
$\phi$ near $z_1$ (near $z_2$ situation is the same by the symmetry)
$$
\phi^{(k)}=c_{-k}(z-z_1)^{-k+1/2}+c_{-k+1}(z-z_1)^{-k+3/2}+\dots+c_0+\dots\ .
$$
This means that for $s>1$ derivative can be estimated by
$$
|c_{-k}|(s-1)^{-k+1/2}+|c_{k+1}|(s-1)^{-k+3/2}+\dots\ .
$$
The series converges in a disc of
a fixed radius (radius is $|z_1+1|$), so its tail can be estimated by a sum
of a geometric progression. Writing these power series explicitly we find
(for $s=1.002$)
\begin{alignat*}{3}
|\phi'| & <  55, & \quad
|\phi''| & <  11800, & \quad
|\phi^{(3)}| &<  8.69\times 10^6,\\
|\phi^{(4)}| & <  1.08 \times 10^{10},  &\quad
|\phi^{(5)}| & <  1.90 \times 10^{13}, &\quad
|\phi^{(6)}| & <  4.25 \times 10^{16}, \\
|\phi^{(7)}| & <  1.17 \times 10^{20}.  &\quad &  & &
\end{alignat*}

The maximal values for first six derivatives of $\nu$ are
\begin{alignat*}{3}
|\nu'|&<0.28, &\quad
|\nu''|&<0.45, &\quad
|\nu^{(3)}|&<1.12, \\
|\nu^{(4)}|&<3.69, &\quad
|\nu^{(5)}|&<15.3,  &\quad
|\nu^{(6)}|&<76.2.
\end{alignat*}

The derivative $\partial_\theta|\phi|$ can be estimated by $r |\phi'|$. We can write
sixth derivative of $\nu(|\phi|)|\phi'|/|\phi|$ as a rational function of
partial derivatives of $|\phi|$, $|\phi'|$, and $\nu$.
 Than we apply triangle inequality and plug in the above estimates.
Finally we have
$$
\left|\frac{\partial\br{\frac{\nu(|\phi|)|\phi'|}{|\phi|}}}{\partial \theta^6}\right|<1.65\times 10^{21}.
$$
Plugging  in the value $\epsilon=\pi/5000$ and the estimate on the sixth derivative into \eqref{error1}
we find that error term in this case is less than $0.0034$.

Next we have to estimate modulus of continuity with respect to $r$.
First we calculate
$$
\partial_r|z-(a+b i)|^2=2 r - 2 (a\cos\theta+b\sin\theta).
$$
Applying this formula several times we find
$$
\begin{aligned}
\partial_r\br{\frac{|\phi'|}{|\phi|}}=
\partial_r \br{\frac{|z-1|}{r\sqrt{|z-z_1||z-z_2|}}} \le
\partial_r \br{\frac{|z-1|}{\sqrt{|z-z_1||z-z_2|}}}
\\
=
\frac{r-\cos\theta}{|z-1|S}-\frac{r-x\cos\theta-y\sin\theta}{2|z-z_1|^2S}|z-1|-
\frac{r-x\cos\theta+y\sin\theta}{2|z-z_2|^2S}|z-1|,
\end{aligned}
$$
where $S=\sqrt{|z-z_1||z-z_2|}$.
Factoring out
$$
\frac{1}{2|z-1|\cdot|z-z_1|^{5/2}|z-z_2|^{5/2}}
$$ we get

$$
\begin{aligned}
2(r-\cos\theta)|z-z_1|^2|z-z_2|^2
-(r-x\cos\theta-y\sin\theta)|z-1|^2|z-z_2|^2
\\
-(r-x\cos\theta+y\sin\theta)|z-1|^2|z-z_1|^2
\\
=-2(r^2-1)(2\cos^2\theta r (x-1)+\cos\theta (r^2+1)(x-1)+2 r y^2).
\end{aligned}
$$
This is a quadratic function with respect to $\cos\theta$. Taking values of $x$ and $y$ into
account we can write it as
$$
\cos^2\theta+\cos\theta\br{r+\frac{1}{r}}\frac{1}{2}-\frac{592}{5625}~.
$$
This quadratic function has two real roots. Their average is
$-(r+1/r)/2<-1$, hence one root is definitely less than $-1$. The
product of roots is a small negative number, which means that the
second root is  positive and less than $1$. Simple calculation shows
that this root decreases as $r$ grows. This means that the
corresponding value of $\theta$ increases. Hence it attains its
maximal value at $r=1.4$ and the maximal value is at most $1.48$.
This gives us that the radial derivative of $|\phi'|/|\phi|$ can be
positive only on the arc $\theta\in[-1.48,1.48]$. By subharmonicity
it attains the maximal on the boundary of $\{z\mid 1<r<1.4,\
-1.48<\theta<1.48\}$. It is not very difficult to check that maximum
is at $z=1.4$ and it is equal to $0.36$.

Let
$$
I(r)=\int_{-\pi}^\pi \nu(|\phi(r^{1/k}e^{i\theta})|
\left|\frac{\phi'(r^{1/k}e^{i\theta})}{\phi(r^{1/k}e^{i\theta})}\right| \frac{d \theta}{2\pi}.
$$
The derivative is
$$
I'(r)=\frac{1}{k r^{1-1/k}}\br{\int_{-\pi}^\pi
 \nu'(|\phi|)\partial_r |\phi| \frac{|\phi'|}{|\phi|}\frac{d \theta}{2\pi}+
\int_{-\pi}^\pi \nu(|\phi|)\partial_r\br{\frac{|\phi'|}{|\phi|}}\frac{d \theta}{2\pi}}.
$$
By the symmetry the first integral is zero. In the second integral
$$
\nu(|\phi|)\partial_r\br{\frac{|\phi'|}{|\phi|}}
$$
can be positive only when $\theta\in[-1.48,1.48]$ and even in this case it is
bounded by $0.36/(r^{1-1/k} k 2\pi)$. Hence
$$
I'(r)<2 \cdot 1.48 \cdot 0.36/(r^{1-1/k}k 2\pi)<0.0131 r^{1/k-1}.
$$
If we compute values $I(r_1)$ and $I(r_2)$ (with precision $0.0034$)
then the minimum of $P(\nu)/\nu$ on $[r_1,r_2]$ is at least
\begin{equation}
r_1^{1/k}(\min\{I(r_1),I(r_2)\}-0.0034-0.0131(r_2-r_1) r^{1/k-1})/\nu(r_1).
\label{int}
\end{equation}

We take $3000$ equidistributed points on $[1,R]$ and compute $I(r)$
at these points. Applying the error estimate (\ref{int}) we find a
rigorous estimate from below of $P\nu/\nu$. The minimum of
$P\nu/\nu$ is at least $1.8079$ which means that
$$
\beta(1)>0.2308.
$$

\def\cprime{$'$}

%\bibliography{snow}

\begin{thebibliography}{1}

\bibitem{BeSmECM}
D.~Beliaev and S.~Smirnov.
\newblock Harmonic measure on fractal sets.
\newblock {\em Proceedings of the 4th European congress of mathematics}, 2005.

\bibitem{BeSmsnow}
D.~Beliaev and S.~Smirnov.
\newblock Random conformal snowflakes.
\newblock Preprint, 2006.

\bibitem{Brennan}
J.~E. Brennan.
\newblock The integrability of the derivative in conformal mapping.
\newblock {\em J. London Math. Soc. (2)}, 18(2):261--272, 1978.

\bibitem{CaJo}
L.~Carleson and P.~W. Jones.
\newblock On coefficient problems for univalent functions and conformal
  dimension.
\newblock {\em Duke Math. J.}, 66(2):169--206, 1992.

\bibitem{Makarov}
N.~G. Makarov.
\newblock Fine structure of harmonic measure.
\newblock {\em St. Petersburg Math. J.}, 10(2):217--268, 1999.

\bibitem{Pommerenke75}
C.~Pommerenke.
\newblock {\em Univalent functions}.
\newblock Vandenhoeck \& Ruprecht, G\"ottingen, 1975.

\bibitem{Pommerenke92}
C.~Pommerenke.
\newblock {\em Boundary behaviour of conformal maps}, volume 299 of {\em
  Grundlehren der Mathematischen Wissenschaften [Fundamental Principles of
  Mathematical Sciences]}.
\newblock Springer-Verlag, Berlin, 1992.

\bibitem{Pommerenke97}
C.~Pommerenke.
\newblock The integral means spectrum of univalent functions.
\newblock {\em Zap. Nauchn. Sem. S.-Peterburg. Otdel. Mat. Inst. Steklov.
  (POMI)}, 237(Anal. Teor. Chisel i Teor. Funkts. 14):119--128, 229, 1997.

\end{thebibliography}
%\bibliographystyle{abbrv}

\end{document}